\documentclass[leqno,draft]{article}

\title{Converse estimates for the simultaneous approximation by Bernstein polynomials with integer coefficients}
\author{Borislav R. Draganov}
\date{}

\usepackage{amsmath,amsthm,amsfonts}
\usepackage[mathscr]{eucal}
\usepackage[english]{babel}
\usepackage{fancyhdr}

\allowdisplaybreaks[1]

\newtheorem{thm}{Theorem}[section]
\newtheorem{prop}[thm]{Proposition}
\newtheorem{cor}[thm]{Corollary}
\newtheorem{lem}[thm]{Lemma}
\theoremstyle{remark}

\numberwithin{equation}{section}

\newcommand{\thmref}[1]{Theorem~\ref{#1}}
\newcommand{\propref}[1]{Proposition \ref{#1}}

\newcommand{\lemref}[1]{Lemma \ref{#1}}

\newcommand{\N}{\mathbb{N}}

\newcommand{\Z}{\mathbb{Z}}
\newcommand{\la}{\left\langle}
\newcommand{\ra}{\right\rangle}

\begin{document}

\maketitle
\bigskip

\thispagestyle{fancy}
\fancyhf{}
\renewcommand{\headrulewidth}{0pt}
\lhead{}
\renewcommand{\footrulewidth}{0.3pt}
\lfoot{\footnotesize This work was supported by grant DN 02/14 of the Fund for Scientific Research of the Bulgarian Ministry of Education and Science.}

\begin{abstract}
We prove a weak converse estimate for the simultaneous approximation by several forms of the Bernstein polynomials with integer coefficients. It is stated in terms of moduli of smoothness. In particular, it yields a big $O$-characterization of the rate of that approximation. We also show that the approximation process generated by these Bernstein polynomials with integer coefficients is saturated. We identify its saturation rate and the trivial class.
\end{abstract}

\bigskip
\noindent
{\footnotesize \leftskip25pt \rightskip25pt {\sl  AMS} {\it classification}: 41A10, 41A25, 41A27, 41A28, 41A29, 41A35, 41A36, 41A40.\\[2pt]
{\it Key words and phrases}: Bernstein polynomials, integer coefficients, integral coefficients, simultaneous approximation, rate of convergence, converse estimates, saturation, modulus of smoothness.
\par}
\bigskip

\section{Main results}

The Bernstein polynomials are defined for $f\in C[0,1]$, $x\in [0,1]$ and $n\in\N_+$ by
\[
B_n f(x):=\sum_{k=0}^n f\left(\frac{k}{n}\right)p_{n,k}(x),\quad p_{n,k}(x):=\binom{n}{k}x^k (1-x)^{n-k}.
\]
It is known that if $f\in C[0,1]$, then
\[
\lim_{n\to\infty} \|B_n f-f\|=0,
\]
where $\|\circ\|$ is the sup-norm on the interval $[0,1]$. The rate of this convergence can be estimated by the Ditzian-Totik modulus of smoothness $\omega_\varphi^2(f,t)$ of the second order with a varying step, controlled by the weight $\varphi(x):=\sqrt{x(1-x)}$, in the uniform norm on the interval $[0,1]$. This modulus is defined by (see \cite[Chapter 2]{Di-To:Mod})
\[
\omega_\varphi^2(f,t):=\sup_{0<h\le t} \|\bar{\Delta}^2_{h\varphi} f\|,
\]
where
\[
\bar{\Delta}^2_{h\varphi(x)} f(x):=\begin{cases}
f(x+h\varphi(x))-2f(x)+f(x-h\varphi(x)), &x\pm h\varphi(x)\in [0,1],\\
0, &\text{otherwise}.
\end{cases}
\]

It was shown that for all $f\in C[0,1]$ and $n\in\N_+$ there holds 
(see \cite{Kn-Zh} and \cite{To}, or \cite[Chapter 10, (7.3)]{De-Lo:CA}, or \cite[Theorem 6.1]{Bu:BP})
\begin{equation}\label{eqBequiv}
c^{-1}\omega_\varphi^2(f,n^{-1/2})\le\|B_n f-f\|\le c\,\omega_\varphi^2(f,n^{-1/2}).
\end{equation}
Throughout $c$ denotes positive constants, whose value is independent of $f$ and $n$. Instead of $\omega_\varphi^2(f,t)$ we can use the moduli defined and considered in \cite{Iv:Dir,Iv:Char}, \cite{Dz-Sh,Ko-Le-Sh:1,Ko-Le-Sh:2,Ko-Le-Sh:3,Ko-Le-Sh:4,Ko-Le-Sh:5,Sh}, or \cite{Dr-Iv:Char}.

Being a linear positive polynomial operator, $B_n$ cannot approximate a function too fast, no matter how ``good'' the function is. Moreover, $B_n$ possesses the property of saturation. More precisely, as \eqref{eqBequiv} and the properties of the modulus $\omega_\varphi^2(f,t)$ show, $\|B_n f-f\|$ cannot tend to $0$ faster than $1/n$ except if $f$ is a linear function, in which case we have $B_n f=f$ for all $n$. Thus the saturation rate of the Bernstein operator is $1/n$, its saturation class consists of those continuous functions $f$ such that $\omega_\varphi^2(f,t)=O(t^2)$, and its trivial class is the set of the linear functions. Let us recall that, by virtue of \cite[Theorem 4.2.1(b)]{Di-To:Mod}, we have for $f\in C[0,1]$
\begin{equation}\label{modchar}
\omega_\varphi^2(f,t)=O(t^2)\ \Longleftrightarrow\ f\in AC[0,1],\ f'\in AC_{loc}(0,1),\ \varphi^2f''\in L_\infty[0,1].
\end{equation}

As is known, the Bernstein operator possesses the property of simultaneous approximation. This means that, if $f\in C^s[0,1]$, $s\in\N_+$, then not only $\|B_n f-f\|\to 0$ as $n\to\infty$, but also $\|(B_n f)^{(i)}-f^{(i)}\|\to 0$, $i=1,\dotsc,s$ (see e.g.~\cite[Chapter 10, Theorem 2.1]{De-Lo:CA}). The rate of this convergence was characterized in \cite{Dr}. In particular, Theorems 1.1 and 1.3 there with $p=\infty$ and $r=1$ imply that the approximation process $(B_n f)^{(s)}\to f^{(s)}$ in uniform norm as $n\to\infty$ is saturated with the rate $1/n$, the trivial class is the set of the algebraic polynomials of degree at most $\max\{1,s-1\}$, and the saturation class consists of the functions $f\in C^s[0,1]$ such that
\[
\omega_\varphi^2(f^{(s)},t)=O(t^2)\quad\text{and}\quad \omega_1(f^{(s)},t)=O(t),
\]
where
\[
\omega_1(F,t):=\sup_{\substack{|x-y|\le t\\x,y\in [0,1]}}|F(x)-F(y)|
\]
is the usual modulus of continuity in the uniform norm on the interval $[0,1]$.

In the present paper we will extend partially the above results to several forms of the Bernstein polynomials with integer coefficients. 

Kantorovich \cite{Ka} (or e.g.\ \cite[pp.\ 3--4]{Fe}, or \cite[Chapter 2, Theorem 4.1]{Lo-Go-Ma:CA}) first introduced such a modification of $B_n$. He considered the operator
\[
\widetilde{B}_n(f)(x):=\sum_{k=0}^n \left[f\left(\frac{k}{n}\right) \binom{n}{k}\right] x^k (1-x)^{n-k}.
\]
Above $[\alpha]$ denotes the largest integer that is less than or equal to the real $\alpha$. 

In \cite{Dr:SimBernInt} we considered another integer form of $B_n$. It is given by
\[
\widehat{B}_n(f)(x):=\sum_{k=0}^n \la f\left(\frac{k}{n}\right) \binom{n}{k} \ra x^k (1-x)^{n-k},
\]
where $\la\alpha\ra$ denotes the nearest integer to the real $\alpha$. More precisely, if $\alpha\ne m+1/2$, $m\in\Z$, we set $\la\alpha\ra$ to be the integer at which $\min_{m\in\Z} |\alpha-m|$ is attained. If $\alpha=m+1/2$, $m\in\Z$, we set either $\la\alpha\ra:=m$, or $\la\alpha\ra:=m+1$ as the definition may depend on whether $m$ is positive or negative, even or odd. The results we will prove are valid regardless of our choice in this case.

We write $\widetilde{B}_n(f)$ and $\widehat{B}_n(f)$, rather than $\widetilde{B}_n f$ and $\widehat{B}_n f$, in order to emphasize that these operators are not linear.

Kantorovich \cite{Ka} showed that, if $f\in C[0,1]$ and $f(0),f(1)\in\Z$, then
\[
\|\widetilde{B}_n(f)-B_n f\|\le\frac{1}{n}.
\]
Similarly, we have 
\begin{equation*}
\|\widehat{B}_n(f)-B_n f\|\le\frac{1}{2n}.
\end{equation*}
Now, applying \eqref{eqBequiv}, we arrive at the characterization
\begin{align}
	&c^{-1}\left(\omega_\varphi^2(f,n^{-1/2})+\frac{1}{n}\right)\le \|\widetilde{B}_n(f)-f\|+\frac{1}{n}
	\le c\left(\omega_\varphi^2(f,n^{-1/2})+\frac{1}{n}\right)\label{charBt}
	\intertext{and}
	&c^{-1}\left(\omega_\varphi^2(f,n^{-1/2})+\frac{1}{n}\right)\le \|\widehat{B}_n(f)-f\|+\frac{1}{n}
	\le c\left(\omega_\varphi^2(f,n^{-1/2})+\frac{1}{n}\right)\label{charBh}
\end{align}
valid for all $f\in C[0,1]$ with $f(0),f(1)\in\Z$.

Consequently, if $0<\alpha\le 1$, then
\begin{align}
	&\|\widetilde{B}_n(f)-f\|=O(n^{-\alpha})\quad\Longleftrightarrow\quad
	\omega_\varphi^2(f,h)=O(h^{2\alpha})\notag
	\intertext{and}
	&\|\widehat{B}_n(f)-f\|=O(n^{-\alpha})\quad\Longleftrightarrow\quad
	\omega_\varphi^2(f,h)=O(h^{2\alpha}),\label{eqhBa}
\end{align}
provided that $f\in C[0,1]$ and $f(0),f(1)\in\Z$. Moreover, as we will prove in \thmref{thmtriv} below, the approximation generated by $\widetilde{B}_n$ and $\widehat{B}_n$ is saturated with the saturation rate of $1/n$ and if $\|\widetilde{B}_n(f)-f\|=o(1/n)$ or $\|\widehat{B}_n(f)-f\|=o(1/n)$, then, similarly to the Bernstein operator, we have that $\widetilde{B}_n(f)=\widehat{B}_n(f)=f$ and $f$ is a polynomial of the type $px+q$, where $p,q\in\Z$.

Here we will also establish analogues of these results for the simultaneous approximation by the operators $\widetilde{B}_n$ and $\widehat{B}_n$.

In \cite{Dr:SimBernInt} we proved direct inequalities for the simultaneous approximation by $\widetilde{B}_n$ and $\widehat{B}_n$. Here we will complement them with the following weak converse estimate.

\begin{thm}\label{thmconv}
Let $s\in\mathbb{N}_+$ and $0<\alpha<1$. Let $f\in C^s[0,1]$, $f(0),f(1)\in\Z$, and
\[
\|(\widetilde{B}_n(f))^{(s)}-f^{(s)}\|=O(n^{-\alpha})\quad\text{or}\quad
\|(\widehat{B}_n(f))^{(s)}-f^{(s)}\|=O(n^{-\alpha}).
\]
Then
\[
\omega_\varphi^2(f^{(s)},h)=O(h^{2\alpha})\quad\text{and}\quad\omega_1(f^{(s)},h)=O(h^{\alpha}).
\]
\end{thm}

Combining this theorem with \cite[Theorems 1.1 and 1.2]{Dr:SimBernInt}, we get the following two big $O$-equivalence relations.

\begin{cor}\label{corequiv1}
Let $s\in\mathbb{N}_+$ and $0<\alpha<1$. Let $f\in C^s[0,1]$ be such that $f(0),f(1),f'(0),f'(1)\in\Z$ 
and $f^{(i)}(0)=f^{(i)}(1)=0$, $i=2,\dots,s$. Let also there exist $n_0\in\N_+$, $n_0\ge s$, such that
\begin{align*}
	f\left(\frac{k}{n}\right) &\ge f(0)+\frac{k}{n}\,f'(0),\quad k=1,\dotsc,s,\ n\ge n_0,\\
	f\left(\frac{k}{n}\right) &\ge f(1)-\left(1-\frac{k}{n}\right)f'(1),\quad k=n-s,\dotsc,n-1,\ n\ge n_0.
\end{align*}
Then
\begin{multline*}
	\|(\widetilde{B}_n(f))^{(s)}-f^{(s)}\|=O(n^{-\alpha})\\
	\Longleftrightarrow\quad
	\omega_\varphi^2(f^{(s)},h)=O(h^{2\alpha})\quad\text{and}\quad\omega_1(f^{(s)},h)=O(h^{\alpha}).
\end{multline*}
\end{cor}

\begin{cor}\label{corequiv2}
Let $s\in\mathbb{N}_+$ and $0<\alpha<1$. Let $f\in C^s[0,1]$ be such that $f(0),f(1),f'(0),f'(1)\in\Z$ 
and $f^{(i)}(0)=f^{(i)}(1)=0$, $i=2,\dots,s$. Then
\begin{multline*}
	\|(\widehat{B}_n(f))^{(s)}-f^{(s)}\|=O(n^{-\alpha})\\
	\Longleftrightarrow\quad
	\omega_\varphi^2(f^{(s)},h)=O(h^{2\alpha})\quad\text{and}\quad\omega_1(f^{(s)},h)=O(h^{\alpha}).
\end{multline*}
\end{cor}

Let us note that the assumptions made in the corollaries are also necessary in order to have simultaneous approximation (see \cite[Theorems 3.1 and 3.2]{Dr:SimBernInt}).

We will also establish the following result, which shows that the approximation processes $(\widetilde{B}_n(f))^{(s)}\to f^{(s)}$ and $(\widehat{B}_n(f))^{(s)}\to f^{(s)}$ in uniform norm are saturated with the saturation rate of $1/n$ and the trivial class consists of the polynomials of the form $px+q$ with $p,q\in\Z$. Note that these processes are neither linear, nor positive.

\begin{thm}\label{thmtriv}
Let $s\in\N_0$ and $f\in C^s[0,1]$ be such that $f(0),f(1)\in\Z$. If
\[
\|(\widetilde{B}_n(f))^{(s)}-f^{(s)}\|=o(1/n)\quad\text{or}\quad
\|(\widehat{B}_n(f))^{(s)}-f^{(s)}\|=o(1/n),
\]
then $f(x)=px+q$ with some $p,q\in\Z$ and thus $\widetilde{B}_n(f)=\widehat{B}_n(f)=f$ for all $n$.
\end{thm}

By virtue of the last theorem with $s=0$, \eqref{charBt}-\eqref{charBh} and \eqref{modchar}, we get the following assertion about the saturation class of the integer forms $\widetilde{B}_n$ and $\widehat{B}_n$ of the Bernstein operator.

\begin{cor}
The operators $\widetilde{B}_n$ and $\widehat{B}_n$ are saturated with the saturation rate of $1/n$. Their saturation class consists of those functions $f\in AC[0,1]$ such that $f(0),f(1)\in\Z$, $f'\in AC_{loc}(0,1)$ and $\varphi^2 f''\in L_\infty[0,1]$.
\end{cor}

I was not able to identify the saturation class of the approximation processes $(\widetilde{B}_n(f))^{(s)}\to f^{(s)}$ and $(\widehat{B}_n(f))^{(s)}\to f^{(s)}$ with $s\ge 1$. In the proof of \thmref{thmtriv} we will note that $(\widetilde{B}_n(f))^{(s)}(x)$ and $(\widetilde{B}_n(f))^{(s)}(x)$ interpolate $f^{(s)}(x)$ at $0$ and $1$ for large $n$, depending on $f$. Therefore the description of the saturation class of these approximation processes might not involve the classical modulus of continuity of $f^{(s)}$ as in Corollaries \ref{corequiv1} and \ref{corequiv2}. However, under an additional assumption, it is quite straightforward to establish the following converse result.

\begin{prop}\label{pr}
Let $s\in\mathbb{N}_+$. Let $f\in C^s[0,1]$, $f(0),f(1)\in\Z$, and $f^{(s)}(x)$ is absolutely continuous with an essentially bounded derivative in some neighbourhoods of $0$ and $1$. If
\[
\|(\widetilde{B}_n(f))^{(s)}-f^{(s)}\|=O(n^{-1})\quad\text{or}\quad
\|(\widehat{B}_n(f))^{(s)}-f^{(s)}\|=O(n^{-1}),
\]
then
\[
\omega_\varphi^2(f^{(s)},h)=O(h^2)\quad\text{and}\quad\omega_1(f^{(s)},h)=O(h);
\]
hence $f^{(s)}\in AC[0,1]$, $f^{(s+1)}\in AC_{loc}(0,1)$ and $f^{(s+1)},\varphi^2 f^{(s+2)}\in L_\infty[0,1]$.
\end{prop}

The contents of the paper are organized as follows. In the next section we will establish the converse estimates formulated in \thmref{thmconv}. The third and last section contains the proofs of \thmref{thmtriv} and \propref{pr}.

\section{Converse estimates}

We will make use of the relation between each of the operators $\widetilde B_n$ and $\widehat B_n$ with $B_n$. 
In \cite[Theorems 2.1 and 2.3]{Dr:SimBernInt} we showed that under the assumptions in Corollaries \ref{corequiv1} and \ref{corequiv2} we have respectively
\begin{align}
	&\|(B_n f)^{(s)}-(\widetilde{B}_n(f))^{(s)}\|\le c\left(\omega_1( f^{(s)},n^{-1}) + \frac{1}{n}\right),\quad n\ge n_0,\label{eqtB-B}
	\intertext{and}
	&\|(B_n f)^{(s)}-(\widehat{B}_n(f))^{(s)}\|\le c\left(\omega_1( f^{(s)},n^{-1}) + \frac{1}{n}\right),\quad n\ge 1.\label{eqhB-B}
\end{align}

Let $s\in\N_+$ and $f\in C^s[0,1]$. Theorems 1.1 and 1.3 in \cite{Dr} with $r=1$ and $p=\infty$, in view of \cite[Theorem 2.1.1]{Di-To:Mod}, imply the strong converse inequalities
\begin{align}
		&\omega_\varphi^{2}(f^{(s)},n^{-1/2})\le c\left(\|(B_{n} f)^{(s)} - f^{(s)}\|+\|(B_{Rn} f)^{(s)} - f^{(s)}\|\right)\label{thmconvsim2}
		\intertext{and}
		&\omega_1(f^{(s)},n^{-1})\le c\left(\|(B_{n} f)^{(s)} - f^{(s)}\|+\|(B_{Rn} f)^{(s)} - f^{(s)}\|\right)\label{thmconvsim1}
\end{align}
for $n\ge n_0$ with some positive integers $R$ and $n_0$, which are independent of $f$ and $n$. It was shown in \cite[Theorem 1.1]{Dr:SimConv} that the two estimates above still hold true without the second term on the right-hand side for $s\le 6$.

Next, we introduce several notations. We will denote the supremum norm of $F$ on the interval $J$ by $\|F\|_J$. When $J=[0,1]$, we will just write $\|F\|$. We set
\begin{align*}
	\tilde b_n(k)&:=\tilde b_n^f(k):=\left[f\left(\frac{k}{n}\right)\binom{n}{k} \right]\,\binom{n}{k}^{-1}
	\intertext{and}
	\hat b_n(k)&:=\hat b_n^f(k):=\la f\left(\frac{k}{n}\right)\binom{n}{k} \ra\,\binom{n}{k}^{-1},
\end{align*} 
where $k=0,\dotsc,n$. Then the operators $\widetilde{B}_n$ and $\widehat{B}_n$ can be written respectively in the form
\begin{align*}
	\widetilde{B}_n (f)(x) &=\sum_{k=0}^n \tilde b_n(k)\,p_{n,k}(x)
	\intertext{and}
	\widehat{B}_n (f)(x) &=\sum_{k=0}^n \hat b_n(k)\,p_{n,k}(x).
\end{align*}

We will use the forward finite difference operator $\Delta_h$ with step $h$, defined by
\[
\Delta_h f(x):=f(x+h)-f(x),\quad \Delta_h^s:=\Delta_h(\Delta_h^{s-1}).
\]
The expanded form of $\Delta_h^s$ is
\begin{equation*}
\Delta_h^s f(x)=\sum_{i=0}^s (-1)^i \binom{s}{i} f(x+(s-i)h),\quad x\in [0,1-sh].
\end{equation*}
We also put $\Delta:=\Delta_1$. Thus we have
\begin{equation*}
\Delta^s \tilde b_n(k)=\sum_{i=0}^s (-1)^i \binom{s}{i} \tilde b_n(k+s-i),\quad k=0,\dotsc,n-s;
\end{equation*}
and analogously for $\hat b_n$.

Let $s\in\N_+$ and $n\ge s$. As is known, the derivatives of $B_n f$ are given by the formula (see \cite{Ma}, or \cite[Chapter 10, (2.3)]{De-Lo:CA})
\begin{equation}\label{eq1}
(B_n f)^{(s)}(x)=\frac{n!}{(n-s)!} \sum_{k=0}^{n-s} \Delta_{1/n}^{s} f\left(\frac{k}{n}\right) p_{n-s,k}(x), \quad x\in [0,1].
\end{equation}

Similarly, we have
\begin{align}
(\widetilde{B}_n (f))^{(s)}(x)&=\frac{n!}{(n-s)!} \sum_{k=0}^{n-s} \Delta^{s} \tilde b_n(k)\,p_{n-s,k}(x), \quad x\in [0,1],\label{eq2a}
\intertext{and}
(\widehat{B}_n (f))^{(s)}(x)&=\frac{n!}{(n-s)!} \sum_{k=0}^{n-s} \Delta^{s} \hat b_n(k)\,p_{n-s,k}(x), \quad x\in [0,1].\label{eq2}
\end{align}

The operators $\widehat B_n$ and $\widetilde B_n$ are not linear. We will use the following property to compensate that. It also incorporates a Bernstein-type inequality.

\begin{lem}\label{lemlinB}
Let $s\in\mathbb{N}_+$, $f\in C^s[0,1]$ and $g\in C^{s+1}[0,1]$. Let $f(0),f(1),\linebreak[1]f'(0),f'(1)\in\Z$ and 
$f^{(i)}(0)=f^{(i)}(1)=0$, $i=2,\dots,s$. Then
\[
\|(\widehat{B}_n(f))^{(s+1)}-(B_n g)^{(s+1)}\|\le c\,n\left(\|f^{(s)}-g^{(s)}\|+\frac{1}{n}\,\|g^{(s+1)}\|+\frac{1}{n}\right),\ n\in\N.
\]
If also there exists $n_0\in\N_+$, $n_0\ge s$, such that for $n\ge n_0$ there hold
\begin{align*}
	&f\left(\frac{k}{n}\right) \ge f(0)+\frac{k}{n}\,f'(0),\quad k=1,\dotsc,s,\\
	&f\left(\frac{k}{n}\right) \ge f(1)-\left(1-\frac{k}{n}\right)f'(1),\quad k=n-s,\dotsc,n-1,
\end{align*}
then
\[
\|(\widetilde{B}_n(f))^{(s+1)}-(B_n g)^{(s+1)}\|\le c\,n\left(\|f^{(s)}-g^{(s)}\|+\frac{1}{n}\,\|g^{(s+1)}\|+\frac{1}{n}\right),\  n\ge n_0.
\]
The constant $c$ is independent of $f$, $g$, and $n$.
\end{lem}

\begin{proof}
We will consider in detail only the operator $\widehat B_n$ and indicate, in due course, the minor changes for $\widetilde B_n$.

We assume that $n\ge s+1$ since otherwise the assertion is trivial. We apply \eqref{eq1} and \eqref{eq2} (or \eqref{eq2a} for $\widetilde B_n$) with $s+1$ in place of $s$, and the identities $\sum_{j=0}^{s+1} \binom{s+1}{j}=2^{s+1}$ and $\sum_{k=0}^{n-s-1} p_{n-s-1,k}(x)\equiv 1$ to deduce for $x\in [0,1]$ that
\begin{equation*}
\begin{split}
&|(\widehat{B}_n(f))^{(s+1)}(x)-(B_n g)^{(s+1)}(x)|\\
&\qquad\le n^{s+1}\sum_{k=0}^{n-s-1} \left|\Delta^{s+1} \hat{b}_n^f(k)-\Delta_{1/n}^{s+1}g\left(\frac{k}{n}\right)\right| p_{n-s-1,k}(x)\\
&\qquad\le n^{s+1}\sum_{k=0}^{n-s-1} \left|\Delta^{s+1} \hat{b}_n^f(k)-\Delta_{1/n}^{s+1}f\left(\frac{k}{n}\right)\right| p_{n-s-1,k}(x)\\
&\qquad\qquad + n^{s+1}\sum_{k=0}^{n-s-1} \left|\Delta_{1/n}^{s+1}(f-g)\left(\frac{k}{n}\right)\right| p_{n-s-1,k}(x)\\
&\qquad\le (2n)^{s+1} \max_{k=0,\dotsc,n} \left|f\left(\frac{k}{n}\right)-\hat{b}_n^f(k)\right|+n^{s+1}\|\Delta_{1/n}^{s+1}(f-g)\|_{[0,1-(s+1)\!/n]}.
\end{split}
\end{equation*}
By virtue of \cite[(2.17), (2.18) and (2.22)]{Dr:SimBernInt} (for $\widetilde B_n$ we use \cite[(2.9), (2.10) and (2.15)]{Dr:SimBernInt} instead) and basic properties of the modulus of continuity, we arrive at
\begin{equation*}
\begin{split}
\left|f\left(\frac{k}{n}\right) -\hat b_n^f(k)\right|\
&\le \frac{c}{n^s}\left(\omega_1( f^{(s)},n^{-1})+\frac{1}{n}\right)\\
&\le \frac{c}{n^s}\left(\omega_1( f^{(s)}-g^{(s)},n^{-1})+\omega_1(g^{(s)},n^{-1})+\frac{1}{n}\right)\\
&\le\frac{c}{n^s}\left(\|f^{(s)}-g^{(s)}\|+\frac{1}{n}\,\|g^{(s+1)}\|+\frac{1}{n}\right),\quad k=0,\dotsc,n.
\end{split}
\end{equation*}

To complete the proof it remains to recall that (see e.g.~\cite[p.~45]{De-Lo:CA})
\begin{equation*}
\|\Delta_{1/n}^{s+1} (f-g)\|_{[0,1-(s+1)/n]}\le 2\,\|\Delta_{1/n}^s (f-g)\|_{[0,1-s/n]}\le \frac{2}{n^s}\,\|f^{(s)}-g^{(s)}\|.
\end{equation*}
\end{proof}

Now, we are ready to give the proof of the weak converse estimate.

\begin{proof}[Proof of \thmref{thmconv}]
We will consider in detail only the operator $\widehat B_n$. Just the same arguments, but based on the corresponding properties of $\widetilde B_n$, yield the assertion for it.

Let $\|(\widehat B_n(f))^{(s)} - f^{(s)}\|\le C_f\,n^{-\alpha}$ for $n\ge n_f$ with some constants $C_f>0$ and $n_f\in\N$ that may depend on 
$f$. Henceforward we will denote by $C_f$ positive constants, which may depend on $f$, but not on $n$ and $h$, $\delta$, and $g$ to be specified below.

We have $\lim_{n\to\infty}\|(\widehat B_n(f))^{(s)} - f^{(s)}\|=0$. Since $f(0),f(1)\in\Z$, we have 
$\lim_{n\to\infty}\|\widehat B_n(f) - f\|=0$ too. Now, \cite[Theorem 3.1]{Dr:SimBernInt} implies that $f^{(i)}(0)=f^{(i)}(1)=0$, $i=2,\dots,s$. For $\widetilde B_n$ we apply \cite[Theorem 3.2]{Dr:SimBernInt} instead. Note also that for both operators we have $f'(0),f'(1)\in\Z$ (see \cite[Section 3]{Dr:SimBernInt}). 

Then \eqref{eqhB-B} (or \eqref{eqtB-B} for $\widetilde B_n$), \eqref{thmconvsim2} and the monotonicity of the modulus of continuity on its second argument imply
\begin{align*}
	\omega_\varphi^{2}(f^{(s)},n^{-1/2})&\le c\left(\|(B_{n} f)^{(s)} - f^{(s)}\|+\|(B_{Rn} f)^{(s)} - f^{(s)}\|\right)\\
	&\le c\left(\|(B_{n} f)^{(s)} - (\widehat{B}_n(f))^{(s)}\| + \|(\widehat{B}_n(f))^{(s)}-f^{(s)}\|\right)\\
	&\qquad + c\left(\|(B_{Rn} f)^{(s)} - (\widehat{B}_{Rn}(f))^{(s)}\| + \|(\widehat{B}_{Rn}(f))^{(s)}-f^{(s)}\|\right)\\
	&\le C_f\left(\omega_1( f^{(s)},n^{-1}) + n^{-\alpha}\right).
\end{align*}
Thus, to complete the proof, it suffices to show that
\begin{equation}\label{eqc1}
\omega_1(f^{(s)},h)=O(h^{\alpha})
\end{equation}
and take into account the monotonicity of $\omega_\varphi^{2}(f^{(s)},h)$ on $h$.

We consider the $K$-functional
\[
K(f^{(s)},t):=\inf_{g\in C^{s+1}[0,1]} \{\|f^{(s)}-g^{(s)}\| + t\,\|g^{(s+1)}\|\}.
\]
As is known (see e.g.~\cite[Chapter 6, Theorem 2.4 and its proof]{De-Lo:CA}), 
\[
\omega_1(f^{(s)},t)\le 2K(f^{(s)},t);
\]
hence, to establish \eqref{eqc1}, it is sufficient to show
\begin{equation}\label{eqc6}
K(f^{(s)},h)=O(h^{\alpha}).
\end{equation}

To this end, we will apply a standard argument based on the Berens-Lorentz Lemma (see \cite{Be-Lo}, or e.g.~\cite[Chapter 10, Lemma 5.2]{De-Lo:CA}).

Let $0<h\le\delta\le 1/n_f$. Set $n:=[1/\delta]$. For any $g\in C^{s+1}[0,1]$, we have
\begin{align*}
	K(f^{(s)},h) &\le \|f^{(s)}-(\widehat{B}_n(f))^{(s)}\| + h\,\|(\widehat{B}_n(f))^{(s+1)}\|\\
	&\le C_f\,n^{-\alpha} + h\,\|(\widehat{B}_n(f))^{(s+1)}-(B_n g)^{(s+1)}\| + h\,\|(B_n g)^{(s+1)}\|\\
	&\le C_f\,\delta^\alpha +c\,\frac{h}{\delta}\left(\|f^{(s)}-g^{(s)}\|+\delta\,\|g^{(s+1)}\|+\delta\right),
\end{align*}
where, at the last step, we estimated the second term by \lemref{lemlinB}, and the third by \cite[Proposition 4.1]{Dr} with $s+1$ in place of $s$, $w=1$ and $p=\infty$. The constant $c$ above is independent of $f$, $g$, $h$, and $\delta$, and $C_f$ is a positive constant, which may depend on $f$, but not on $g$, $h$, and $\delta$.

We take the infimum on $g\in C^{s+1}[0,1]$ and thus arrive at
\begin{equation*}
K(f^{(s)},h) + h\le C_f\,\delta^\alpha +c\,\frac{h}{\delta}\left(K(f^{(s)},\delta) + \delta\right).
\end{equation*}
Now, the Berens-Lorentz Lemma with $\phi(x):=K(f^{(s)},x^2) + x^2$ and $2\alpha$ in place of $\alpha$ (in the notations of 
\cite[Chapter 10, Lemma 5.2]{De-Lo:CA}) implies \eqref{eqc6}.
\end{proof}

\section{Saturation}

In this section we will first prove \thmref{thmtriv}. It shows that the approximation processes $(\widetilde{B}_n (f))^{(s)}\to f^{(s)}$ and $(\widehat{B}_n (f))^{(s)}\to f^{(s)}$ are saturated.

\begin{proof}[Proof of \thmref{thmtriv}]
We consider $\widehat{B}_n$. The argument for $\widetilde{B}_n$ is just the same.

First of all, let us note that if $f(x)=px+q$ with $p,q\in\Z$, then 
\[
\left (p\,\frac{k}{n}+q\right)\binom{n}{k}\in\Z,\quad k=0,\dotsc,n;
\]
hence $\widehat{B}_n(f)=B_n f$. As is known, $B_n$ preserves the linear functions. Therefore $\widehat{B}_n(f)=f$ for all $n$.

We consider the case $s=0$. Let $\delta\in (0,1/2)$ be fixed. For $x\in [\delta,1-\delta]$ we have
\begin{align*}
	|B_n f(x)-\widehat{B}_n(f)(x)|
	&\le \sum_{k=1}^{n-1}\left|f\left(\frac{k}{n}\right) \binom{n}{k}-\la f\left(\frac{k}{n}\right) \binom{n}{k} \ra \right|x^k (1-x)^{n-k}\\
	&\le \frac{1}{2}\sum_{k=1}^{n-1} x^k (1-x)^{n-k}\le \frac{1}{2}\sum_{k=1}^{n-1} (1-\delta)^k (1-\delta)^{n-k}\\
	&=\frac{n-1}{2}\,(1-\delta)^n.
\end{align*}
Consequently,
\begin{equation}\label{eqBo}
\|B_n f-f\|_{[\delta,1-\delta]}=o(1/n).
\end{equation}

Further, by virtue of \eqref{eqhBa} with $\alpha=1$ and $\|\widehat{B}_n(f)-f\|=o(1/n)$, we get $\omega_\varphi^2(f,h)=O(h^2)$. Therefore $f\in W^2_\infty[\delta,1-\delta]$ (see \eqref{modchar}).

Now, Voronovskaya's classical result (see e.g.~\cite[Chapter 10, Theorem 3.1]{De-Lo:CA}) and \eqref{eqBo} yield that $f''(x)=0$ a.e.~in $[\delta,1-\delta]$. Since $\delta$ was arbitrarily fixed in $(0,1/2)$, we arrive at $f''(x)=0$ a.e.~in $[0,1]$. Consequently, $f(x)$ is a linear function. It assumes integral values at $0$ and $1$; hence $f(x)=px+q$ with some $p,q\in\Z$.

Let $s\in\N_+$. As is known, for any $g\in C^s[0,1]$ we have (see e.g.~\cite[Chapter 2, Theorem 5.6]{De-Lo:CA})
\[
\|g^{(i)}\|\le c\big(\|g\| + \|g^{(s)}\|\big),\quad i=1,\dotsc,s-1.
\]
Therefore
\[
\lim_{n\to\infty}\|\widehat{B}_n(f)-f\|=0\quad\text{and}\quad
\lim_{n\to\infty}\|(\widehat{B}_n(f))^{(s)}-f^{(s)}\|=0
\]
imply
\[
\lim_{n\to\infty}\|(\widehat{B}_n(f))^{(i)}-f^{(i)}\|=0,\quad i=1,\dotsc,s-1.
\]
In particular, we have $\lim_{n\to\infty}(\widehat{B}_n(f))^{(i)}(0)=f^{(i)}(0)$, $i=0,\dotsc,s-1$. Since $(\widehat{B}_n(f))^{(i)}(0)\in\Z$, we deduce that for all $n$ large enough we have $(\widehat{B}_n(f))^{(i)}(0)=f^{(i)}(0)$, $i=0,\dotsc,s-1$.

Consequently,
\[
\widehat{B}_n(f)(x)-f(x)=\frac{1}{(s-1)!}\int_0^x (x-u)^{s-1}\left((\widehat{B}_n(f))^{(s)}(u)-f^{(s)}(u)\right)du;
\]
hence
\[
\|\widehat{B}_n(f)-f\|=o(1/n),
\]
which reduces the assertion to the case $s=0$.
\end{proof}

\begin{proof}[Proof of \propref{pr}]
We will consider only the operator $\widehat{B}_n$. The proof for $\widetilde{B}_n$ is quite similar.

As in the proof of \thmref{thmconv} we first deduce that $f'(0),f'(1)\in\Z$ and $f^{(i)}(0)=f^{(i)}(1)=0$, $i=2,\dots,s$. Then we observe that the considerations in the proof of \cite[Theorem 2.3]{Dr:SimBernInt} actually imply
\[
\|(B_n f)^{(s)}-(\widehat{B}_n(f))^{(s)}\|\le c\left(\omega_1( f^{(s)},n^{-1})_{[0,s/n]} +\omega_1( f^{(s)},n^{-1})_{[1-s/n,1]}+ \frac{1}{n}\right),
\]
where we have set for the interval $J\subset [0,1]$
\[
\omega_1(F,t)_J:=\sup_{\substack{|x-y|\le t\\x,y\in J}}|F(x)-F(y)|.
\]
We have $f^{(s)}\in W^1_\infty[0,s/n]$ and $f^{(s)}\in W^1_\infty[1-s/n,1]$ for all $n$ large enough; hence
\[
\|(B_n f)^{(s)}-(\widehat{B}_n(f))^{(s)}\|=O(n^{-1}).
\] 

Consequently,
\[
\|(B_n f)^{(s)}-f^{(s)}\|=O(n^{-1}).
\] 
By virtue of \eqref{thmconvsim2}-\eqref{thmconvsim1}, this implies 
\[
\omega_\varphi^2(f^{(s)},t)=O(t^2)\quad\text{and}\quad \omega_1(f^{(s)},t)=O(t).
\]

Basic properties of the moduli (see \eqref{modchar} and \cite[Chapter 2, Theorem 9.3]{De-Lo:CA}) yield the second assertion of the proposition. 
\end{proof}

\bigskip
\begin{footnotesize}
\noindent
\begin{tabular}{ll}
Borislav R. Draganov& \\
Dept. of Mathematics and Informatics&
Inst. of Mathematics and Informatics\\
Sofia University ``St. Kliment Ohridski''&
Bulgarian Academy of Sciences\\
5 James Bourchier Blvd.&
bl. 8 Acad. G. Bonchev Str.\\
1164 Sofia&
1113 Sofia\\
Bulgaria&
Bulgaria\\
bdraganov@fmi.uni-sofia.bg& \\
\end{tabular}

\end{footnotesize}


\begin{thebibliography}{50}

\bibitem{Fe}
Le Baron O. Ferguson, Approximation by Polynomials with Integral Coefficients, Mathematical Surveys Vol. 17, American Mathematical Society, 1980.

\bibitem{Be-Lo}
H. Berens, G. G. Lorentz, Inverse theorems for Bernstein polynomials, Indiana Univ. Math. J. 21 (1972), 693--708.

\bibitem{Bu:BP}
J. Bustamante, Bernstein Operators and Their Properties, Birkh\"auser, 2017.

\bibitem{De-Lo:CA}
R. A. DeVore, G. G. Lorentz, Constructive Approximation, Springer-Verlag, Berlin, 1993.

\bibitem{Di-To:Mod}
Z. Ditzian, V. Totik, Moduli of Smoothness, Springer-Verlag, New York, 1987.

\bibitem{Dr}
B. R. Draganov, Strong estimates of the weighted simultaneous approximation by the Bernstein and Kantorovich operators and their iterated Boolean sums, J. Approx. Theory 200 (2015), 92--135.

\bibitem{Dr:SimConv}
B. R. Draganov, An exact strong converse inequality for the weighted simultaneous approximation by the Bernstein operator, in ``Constructive Theory of Functions, Sozopol 2016'' (K.~Ivanov, G.~Nikolov, R.~Uluchev, Eds.), pp.~75--97, Marin Drinov Academic Publishing House, Sofia, 2018.

\bibitem{Dr:SimBernInt}
B. R. Draganov, Simultaneous approximation by Bernstein polynomials with integer coefficients, J. Approx. Theory 237 (2019), 1--16.

\bibitem{Dr-Iv:Char}
B. R. Draganov, K. G. Ivanov, A New characterization of weighted Peetre $K$-functionals, Constr. Approx. 21 (2005), 113--148.

\bibitem{Dz-Sh}
V. K. Dzyadyk, I. A. Shevchuk, Theory of Uniform Approximation of Functions by Polynomials, Walter de Gruyter, Berlin, 2008.

\bibitem{Iv:Dir}
K. G. Ivanov, Some characterizations of the best algebraic approximation in $L_p[-1,1]$ $(1\le p \le \infty)$, C. R. Acad. Bulgare Sci.
34 (1981), 1229--1232.

\bibitem{Iv:Char}
K. G. Ivanov, A characterization of weighted Peetre $K$-functionals, J. Approx. Theory 56 (1989), 185--211.

\bibitem{Ka}
L. V. Kantorovich, Some remarks on the approximation of functions by means of polynomials with integer coefficients, Izv. Akad. Nauk SSSR, Ser. Mat. 9 (1931), 1163--1168 (in Russian).

\bibitem{Kn-Zh}
H.-B. Knoop and X.-L. Zhou, The lower estimate for linear positive operators (II), Results Math. 25 (1994), 315--330.

\bibitem{Ko-Le-Sh:1}
K. Kopotun, D. Leviatan and I. A. Shevchuk, New moduli of smoothness, Publ. Inst. Math. Serbian Acad. Sci. Arts Belgr. \textbf{96} (110) (2014), 169--180.

\bibitem{Ko-Le-Sh:2}
K. Kopotun, D. Leviatan and I. A. Shevchuk, New moduli of smoothness: weighted DT moduli revisited and applied, Constr. Approx. \textbf{42} (2015), 129--159.

\bibitem{Ko-Le-Sh:3}
K. Kopotun, D. Leviatan and I. A. Shevchuk, On the moduli of smoothness with Jacobi weights, Ukrainian Math. J. \textbf{70} (3) (2018), 437--466.

\bibitem{Ko-Le-Sh:4}
K. Kopotun, D. Leviatan and I.A. Shevchuk, On weighted approximation with Jacobi weights, J. Approx. Theory \textbf{237} (2019), 96--112.

\bibitem{Ko-Le-Sh:5}
K. Kopotun, D. Leviatan and I.A. Shevchuk, On some properties of moduli of smoothness with Jacobi weights, in ``Topics in Classical and Modern Analysis. Applied and Numerical Harmonic Analysis'' (M. Abell, E. Iacob, A. Stokolos, S. Taylor, S. Tikhonov, J. Zhu, Eds.), pp.~19--31, Birkh\"auser, Cham, 2019.

\bibitem{Lo-Go-Ma:CA}
G. G. Lorentz, M. v.Golitschek, Y. Makovoz, Constructive Approximation, Advanced Problems, Springer-Verlag, Berlin, 1996. 

\bibitem{Ma}
R. Martini, On the approximation of functions together with their derivatives by certain linear positive operators, Indag. Math. 31 (1969), 473--481. 

\bibitem{Sh}
I. A. Shevchuk, Polynomial Approximation and Traces of Functions Continuous on a Segment, Naukova Dumka, Kiev, 1992 (in Russian).

\bibitem{To}
V. Totik, Approximation by Bernstein polynomials, Amer. J. Math. 116 (1994), 995--1018.


\end{thebibliography}
\end{document}